\documentclass{article}
\begin{document}
\author{Richard Drociuk}
\title{On the Root Ambiguity in the Complete Solution to the Most General Fifth Degree Polynomial}
\maketitle
\section{Abstract}
Starting from the solution to Bring's equation the root ambiguity is removed from the solution to the quintic equation. This gives the five complex roots of the quintic equation as indicated by Gauss's Fundamental Theorem of Algebra.
\section{Introduction}
In the previous paper[Drociuk,1], the solution to the quintic was given, but the root ambiguity was not correctly removed. This problem arises because Ferrari's method for solving the quartic equation does not introduce all possible roots. Instead if one uses a Tshirnhaussen transformation as demonstrated[Drociuk,1], to the quartic equation, the five correct roots of the quintic equation,\\
\begin{equation}\label{1}
x^{5}+mx^{4}+nx^{3}+px^{2}+qx+r=0
\end{equation}
are selcted for arbitrary coefficients $m$, $n$, $p$, $q$, and $r$. The roots are selected from all possible roots using conditional loops to an arbitrary precission, $\epsilon$.
\section{Ambiguity in the Quartic Solution}
The quartic equation,
\begin{equation}\label{2}
x^{4}+a_{3}x^{3}+a_{2}x^{2}+a_{1}x+a_{0}=0
\end{equation}
is transformed with,
\begin{equation}\label{3}
Tsh2 = x^{3}+b_{2}x^{2}+b_{1}x+b_{0}+y_{n}
\end{equation}
to the quadratic in $y_{n}^{2}$,
\begin{equation}\label{4}
y_{n}^{4}+B_{2}y_{n}^{2}+B_{0}=0
\end{equation}
 whose four roots are,
\begin{equation}\label{5}
y_{1}= {\frac{1}{2}}(-2B_{2}+2(B_{2}^{2}-4B_{0})^{\frac{1}{2}})^{\frac{1}{2}}
\end{equation}
\newpage
\begin{equation}\label{6}
y_{2}= -{\frac{1}{2}}(-2B_{2}+2(B_{2}^{2}-4B_{0})^{\frac{1}{2}})^{\frac{1}{2}}
\end{equation}
\begin{equation}\label{7}
y_{3}= {\frac{1}{2}}(-2B_{2}-2(B_{2}^{2}-4B_{0})^{\frac{1}{2}})^{\frac{1}{2}}
\end{equation}
\begin{equation}\label{8}
y_{4}= -{\frac{1}{2}}(-2B_{2}-2(B_{2}^{2}-4B_{0})^{\frac{1}{2}})^{\frac{1}{2}}
\end{equation}
where,\\ 
\\
\begin{math}
B0 = 3b_{0}^{2}a_{3}^{2}a_{0}+b_{0}^{3}b_{2}a_{3}^{2}+b_{0}^{2}b_{2}^{2}a_{2}^{2}-2b_{0}^{3}b_{2}a_{2}+3a_{3}a_{2}b_{0}^{3}+b_{0}b_{2}^{3}a_{1}^{2}-b_{0}^{3}b_{1}a_{3}
+4b_{0}b_{2}a_{0}^{2}+4b_{0}^{2}b_{1}a_{0}+b_{0}^{2}b_{1}^{2}a_{2}-2b_{0}^{2}b_{1}a_{2}^{2}-b_{2}^{3}a_{0}^{2}a_{3}-2a_{0}b_{1}^{3}a_{2}+a_{3}^{2}a_{0}b_{1}^{3}
-b_{2}a_{0}^{2}a_{1}-2b_{1}a_{0}^{2}a_{2}+b_{1}a_{0}a_{1}^{2}+a_{0}b_{1}^{2}a_{2}^{2}+a_{1}b_{0}^{2}b_{2}a_{2}-3b_{0}a_{3}a_{0}^{2}+2a_{1}b_{0}^{2}b_{2}a_{3}^{2}
-3a_{1}a_{3}a_{2}b_{0}^{2}+b_{0}^{2}a_{3}^{2}a_{2}b_{1}-5b_{0}^{2}b_{2}a_{3}a_{0}-b_{0}^{2}a_{3}a_{2}^{2}b_{2}+3b_{0}^{2}b_{1}b_{2}a_{1}
-2b_{0}^{2}b_{2}^{2}a_{3}a_{1}-2b_{0}b_{2}^{3}a_{0}a_{2}-4b_{0}b_{2}a_{0}b_{1}^{2}-2b_{0}b_{2}a_{0}a_{2}^{2}+2b_{0}a_{1}b_{1}^{2}a_{2}+b_{0}b_{2}^{2}a_{0}a_{1}
-b_{0}a_{3}^{2}a_{1}b_{1}^{2}+2b_{0}a_{3}a_{1}^{2}b_{1}-b_{0}b_{2}^{2}a_{3}a_{1}^{2}+b_{0}b_{2}a_{2}a_{1}^{2}-b_{0}a_{1}b_{1}a_{2}^{2}+3b_{0}a_{1}a_{2}a_{0}
-3b_{0}b_{1}b_{2}a_{1}^{2}-5b_{0}a_{1}b_{1}a_{0}+3b_{1}b_{2}a_{3}a_{0}^{2}+a_{0}b_{0}b_{1}^{2}a_{3}+a_{0}b_{1}a_{3}a_{2}b_{0}-a_{0}b_{1}b_{2}a_{2}a_{1}
-a_{0}b_{2}a_{2}b_{1}^2a_{3}+2b_{1}^{2}a_{0}^{2}+b_{2}^{4}a_{0}^{2}+b_{1}^{4}a_{0}-b_{0}b_{1}^{3}a_{1}+a_{0}^{2}b_{2}^{2}a_{2}-4b_{2}^{2}a_{0}^{2}b_{1}
-3a_{1}b_{0}^{3}+3b_{0}^{2}a_{1}^{2}-b_{0}a_{1}^{3}-3b_{0}^{2}a_{2}a_{0}+b_{0}^{4}+2b_{2}^{2}a_{0}b_{0}^{2}-b_{0}^{3}a_{3}^{3}+b_{0}^{2}a_{2}^{3}-a_{1}b_{0}^{2}b_{1}a_{3}
+a_{0}^{3}-b_{0}^{2}b_{2}a_{2}b_{1}a_{3}-b_{0}a_{1}b_{2}a_{3}a_{0}+b_{0}b_{2}a_{1}b_{1}^{2}a_{3}+2b_{0}a_{0}b_{2}^{2}a_{3}a_{2}-3b_{0}a_{0}b_{1}b_{2}a_{3}^{2}
+4b_{0}b_{2}a_{0}b_{1}a_{2}+3b_{0}b_{2}^{2}a_{0}b_{1}a_{3}-b_{0}b_{1}a_{1}b_{2}^{2}a_{2}+b_{0}a_{1}b_{2}a_{2}b_{1}a_{3}+3b_{1}^{2}a_{0}b_{2}a_{1}
-2b_{1}^{2}a_{0}a_{1}a_{3}-b_{1}a_{0}b_{2}^{3}a_{1}+b_{1}^{2}a_{0}b_{2}^{2}a_{2}+b_{2}^{2}a_{0}b_{1}a_{3}a_{1}-b_{2}a_{0}b_{1}^{3}a_{3}
\end{math}
\begin{equation}\label{9}
\end{equation}
and\\
\\
\begin{math}
B2 = a_{3}^{2}a_{2}b_{1}+2a_{1}b_{2}a_{3}^{2}-6b_{0}b_{2}a_{2}-5b_{2}a_{3}a_{0}+4b_{1}a_{0}-3a_{2}a_{0}+3a_{3}^{2}a_{0}-3a_{3}a_{2}a_{1}-3b_{0}a_{3}^{3}+6b_{0}^{2}-9b_{0}a_{1}+b_{2}^{2}a_{2}^{2}+3a_{1}^{2}-2b_{1}a_{2}^{2}+b_{1}^{2}a_{2}-b_{2}a_{2}b_{1}a_{3}-a_{3}a_{2}^{2}b_{2}+a_{2}^{3}+3b_{1}b_{2}a_{1}-a_{1}b_{1}a_{3}-3b_{0}b_{1}a_{3}+3b_{0}b_{2}a_{3}^{2}+2a_{0}b_{2}^{2}+9a_{3}a_{2}b_{0}+b_{2}a_{2}a_{1}-2b_{2}^{2}a_{3}a_{1}
\end{math}
\begin{equation}\label{10}
\end{equation}
The coefficients of the Tshirnhaussen transformation (\ref{3}) are,
\\
\begin{math}
b_{0} = \frac{(36b_{01}b_{02}b_{03}-108b_{00}b_{03}^{2}-8b_{02}^{3}+12(3)^{\frac{1}{2}}(4b_{01}^{3}b_{03}-b_{01}^{2}b_{02}^{2}-18b_{01}b_{02}b_{03}b_{00}+27b_{00}^{2}b_{03}^{2}+4b_{00}b_{02}^{3})^{\frac{1}{2}}b_{03})^{\frac{1}{3}}}{6b_{03}}
\end{math}
\\
\begin{math}
-{\frac{2}{3}}{\frac{(3b_{01}b_{03}-b_{02}^{2})}{(b_{03}(36b_{01}b_{02}b_{03}-108b_{00}b_{03}^{2}-8b_{02}^{3}+12(3)^{\frac{1}{2}}(4b_{01}^{3}b_{03}-b_{01}^{2}b_{02}^{2}-18b_{01}b_{02}b_{03}b_{00}+27b_{00}^{2}b_{03}^{2}+4b_{00}b_{02}^{3})^{\frac{1}{2}}b_{03})^{\frac{1}{3}}-{\frac{1}{3}}{\frac{b_{02}}{b_{03}}})}}  
\end{math}
\begin{equation}\label{11}
\end{equation}
choose,
\\
\begin{equation}\label{12}
b_{1} = 0
\end{equation}
and
\\
\begin{equation}\label{13}
b_{2} = {\frac{(a_{3}^{3}+3a_{1}-4b_{0}+b_{1}a_{3}-3a_{3}a_{2})}{(a_{3}^{2}-2a_{2})}}
\end{equation}
\newpage
the coefficients $b_{00}$, $b_{01}$, $b_{02}$ and $b_{03}$ are given by\\
\\
\begin{math}
b_{00} = (a_{0}^{2}a_{3}^{7}+20a_{2}^{3}a_{1}^{3}+2a_{3}^{6}a_{1}^{3}+18a_{3}^{3}a_{1}^{4}-36a_{1}a_{2}^{3}a_{0}a_{3}^{2}+150a_{3}a_{1}^{2}a_{2}^{2}a_{0}+29a_{1}a_{2}^{2}a_{3}^{4}a_{0}-54a_{3}^{3}a_{1}^{2}a_{0}a_{2}-4a_{3}^{6}a_{0}a_{1}a_{2}-48a_{3}^{2}a_{0}^{2}a_{1}a_{2}+27a_{1}^{5}-a_{3}^{5}a_{2}^{2}a_{1}^{2}+28a_{3}^{3}a_{0}^{2}a_{2}^{2}-24a_{2}^{3}a_{0}^{2}a_{3}+48a_{2}^{2}a_{0}^{2}a_{1}+24a_{2}^{5}a_{0}a_{3}+12a_{3}^{4}a_{0}^{2}a_{1}-10a_{3}^{5}a_{0}^{2}a_{2}-14a_{3}^{3}a_{0}a_{2}^{4}+2a_{3}^{5}a_{0}a_{2}^{3}-12a_{3}a_{2}^{4}a_{1}^{2}+7a_{3}^{3}a_{2}^{3}a_{1}^{2}-48a_{1}a_{2}^{4}a_{0}+3a_{3}^{5}a_{1}^{2}a_{0}+21a_{2}^{2}a_{1}^{3}a_{3}^{2}-15a_{2}a_{1}^{3}a_{3}^{4}-72a_{1}^{3}a_{2}a_{0}+9a_{1}^{3}a_{3}^{2}a_{0}-63a_{3}a_{1}^{4}a_{2})/(a_{3}^{2}-2a_{2})^{3}
\end{math}
\begin{equation}\label{14}
\end{equation}
\begin{math}
b_{01} = (-10a_{3}^{6}a_{1}^{2}-84a_{3}^{3}a_{1}^{3}+80a_{0}a_{2}^{4}-76a_{2}^{3}a_{1}^{2}-64a_{0}^{2}a_{2}^{2}-2a_{3}^{4}a_{2}^{4}+12a_{3}^{2}a_{2}^{5}-16a_{3}^{4}a_{0}^{2}-108a_{1}^{4}-16a_{2}^{6}+120a_{3}^{3}a_{1}a_{2}a_{0}-344a_{3}a_{1}a_{2}^{2}a_{0}+16a_{2}^{3}a_{3}^{2}a_{0}-10a_{3}^{5}a_{1}a_{0}-24a_{2}^{2}a_{3}^{4}a_{0}-74a_{2}^{2}a_{1}^{2}a_{3}^{2}+70a_{2}a_{1}^{2}a_{3}^{4}+72a_{3}a_{2}^{4}a_{1}-52a_{3}^{3}a_{2}^{3}a_{1}+8a_{3}^{5}a_{2}^{2}a_{1}+4a_{3}^{6}a_{0}a_{2}+64a_{3}^{2}a_{0}^{2}a_{2}240a_{3}a_{1}^{3}a_{2}+12a_{3}^{2}a_{0}a_{1}^{2}+192a_{0}a_{2}a_{1}^{2})/(a_{3}^{2}-2a_{2})^{3}
\end{math}
\begin{equation}\label{15}
\end{equation}
\begin{math}
b_{02} = (-64a_{3}^{3}a_{2}a_{0}+8a_{3}^{5}a_{0}+16a_{1}a_{3}^{6}-104a_{1}a_{3}^{4}a_{2}+112a_{3}^{2}a_{1}a_{2}^{2}-64a_{3}a_{2}^{4}+128a_{1}^{2}a_{3}^{3}-8a_{3}^{5}a_{2}^{2}-304a_{3}a_{2}a_{1}^{2}+144a_{1}^{3}+64a_{1}a_{2}^{3}+48a_{3}^{3}a_{2}^{3}+192a_{0}a_{3}a_{2}^{2}-80a_{1}a_{3}^{2}a_{0}-128a_{1}a_{2}a_{0})/(a_{3}^{2}-2a_{2})^{3}
\end{math}
\begin{equation}\label{16}
\end{equation}
\begin{math}
b_{03} = (-64a_{3}^{2}a_{2}^{2}+48a_{3}^{4}a_{2}-8a_{3}^{6}+64a_{3}^{2}a_{0}-64a_{1}^{2}+128a_{3}a_{2}a_{1}-64a_{1}a_{3}^{3})/(a_{3}^{2}-2a_{2})^{3}
\end{math}
\begin{equation}\label{17}
\end{equation}
Now the cubic(\ref{3})is solved for it's three roots,\\
\\
\begin{math}
x_{1n} = \frac{\Delta_{n}-(2b_{1}-\frac{2}{3}b_{2}^{2})}{6\Delta_{n}-\frac{1}{3}b_{2}}
\end{math}
\begin{equation}\label{18}
\end{equation}
\begin{math}
x_{2n} = -\frac{1}{12}\Delta_{n}+\frac{b_{1}-\frac{1}{3}b_{2}^{2}}{\Delta_{n}}-\frac{1}{3}b_{2}+i\sqrt{3}(\frac{1}{12}\Delta_{n}+\frac{b_{1}-\frac{1}{3}b_{2}^{2}}{\Delta_{n}})
\end{math}
\begin{equation}\label{19}
\end{equation}
\begin{math}
x_{3n} = -\frac{1}{12}\Delta_{n}+\frac{b_{1}-\frac{1}{3}b_{2}^{2}}{\Delta_{n}}-\frac{1}{3}b_{2}-i\sqrt{3}(\frac{1}{12}\Delta_{n}+\frac{b_{1}-\frac{1}{3}b_{2}^{2}}{\Delta_{n}})
\end{math}
\begin{equation}\label{20}
\end{equation}
\begin{math}
\Delta_{n} =(36b_{1}b_{2}-108y_{1}-108b_{0}-8b_{2}^{3}+12(12b_{1}^{3}-3b_{1}^{2}b_{2}^{2}-54b_{1}b_{2}y_{n}-54b_{1}b_{2}b_{0}+81y_{n}^{2}+162y_{n}b_{0}+12y_{n}b_{2}^{3}+81b_{0}^{2}+12b_{0}b_{2}^{3})^{\frac{1}{2}})^{\frac{1}{3}}
\end{math}
\begin{equation}\label{21}
\end{equation}
By substituting the values of $y_{n}$, from equations (\ref{5}) to (\ref{8}), into equations (\ref{18}) to (\ref{20}), one can generate twelve possible roots, $x_{mn}$ with $m \in (1,2,3)$ and $n \in (1,2,3,4)$. Only four roots satisfy equation (\ref{2}), they are found by direct substitution. So for all complex coeffients $a_{0}$, $a_{1}$, $a_{2}$ and $a_{3}$ in equation (\ref{2}), the four roots of the quartic are contained in $x_{mn}$.
\newpage
\section{Ambiguity in the Quintic Solution}
Recall in [Drociuk,1], the coefficients of the quartic Tshirnhausen Transformation,
\begin{equation}\label{22}
Tsh1 = x^{4}+dx^{3}+cx^{2}+bx+a+y
\end{equation} 
and let the coefficients of equation (\ref{2}) be,
\begin{equation}\label{23}
a_{0} = a+y   
\end{equation}
\begin{equation}\label{24}
a_{1} = b   
\end{equation}
\begin{equation}\label{25}
a_{2} = c   
\end{equation}
\begin{equation}\label{26}
a_{3} = d
\end{equation}
and calculate the twelve $x_{mn}$ using the equations of the previous section. Then using nested ``for loops'' and ``if-then'' statements on Maple, $x_{mn}$ is substituted into both the quintic (\ref{1}) and the quartic (\ref{22}). When both these equations are made less than $\epsilon$, we have obtained the first root of the quintc equation (\ref{1}), let it be $r_{1}$. $\epsilon$ can be made arbitrarily close to zero. Let $\epsilon = 10^{-100}$, you may chose it to be zero if you want to wait or have a faster computer. Now $r_{1}$ is then factored out of the quintic (\ref{1}), leaving a quartic equation to be solved with the following coefficients,
\begin{equation}\label{27}
a_{0} = q+r_{1}p+r_{1}^{2}n+mr_{3}+r_{1}^{4}
\end{equation}
\begin{equation}\label{28}
a_{1} = p+r_{1}n+r_{1}^{2}m+r_{1}^{3}
\end{equation}
\begin{equation}\label{29}
a_{2} = n+mr_{1}+r_{1}^{2}
\end{equation}
\begin{equation}\label{30}
a_{3} = m+r_{1}
\end{equation}
Equations (\ref{27}) to (\ref{30}) are then substituted into the solution of the general quartic of the previous section, which gives another twelve roots, $x_{mn}$. This time four of these roots satisfy both the quartic (\ref{2}) with coefficients (\ref{27}) to (\ref{30}) and the quintic (\ref{1}), let them be $r_{2}$,$ r_{3}$, $r_{4}$ and $r_{5}$. They are determined in the same way as $r_{1}$, using nested ``for loops'' and ``if-then'' statements on Maple. The the final array has five non-zero elements, they are the five complex roots of the quintic equation,(\ref{1}).
\section{Conclusion}
The root ambiguity is removed from the solution to the general fifth and fourthdegree polynomials.       
\section{Reference}
1) Drociuk, Richard, ``On the Complete Solution to the Most General Fifth Degree Polynomial'', GM/0005026v1 3 May 2000.
\end{document}